\documentclass{article}
\usepackage{amssymb}
\usepackage{amsmath}

\newcommand \leqo { \mathord{\leq} }
\newcommand \geqo { \mathord{\geq} }
\input{tcilatex}

\begin{document}

\begin{center}
\textbf{G. Cz\'{e}dli's tolerance factor lattice construction }

\textbf{and weak ordered relations}

\medskip

\textsc{S\'{a}ndor Radeleczki\footnote[1]{{\small This research started as
part of the TAMOP-4.2.1.B-10/2/KONV-2010-0001 project, supported by the
European Union, co-financed by the European Social Fund 113/173/0-2}.}}

\bigskip
\end{center}

\textsc{Abstract:} {\small G. Cz\'{e}dli proved that the blocks of any
compatible tolerance }$T${\small \ of a lattice }$L${\small \ can be ordered
in such a way that they form a lattice }$L/T${\small \ called the factor
lattice of }$L${\small \ modulo }$T${\small . Here we show that Dedekind-Mac
Neille completion of the lattice }$L/T${\small \ is isomorphic to the
concept lattice of the context }$(L,L,R)${\small , where }$R${\small \
stands for the reflexive weak ordered relation }$\leqo\circ T\circ \leqo$%
{\small . Weak ordered relations constitute the generalization of the
ordered relations introduced by S. Valentini. Reflexive weak ordered
relations can be characterized as compatible reflexive relations }$%
R\subseteq L^{2}${\small \ satisfying }$R=\ \leqo\circ R\circ \leqo${\small .%
}

\medskip

\textit{2010 Mathematics Subject Classification}. Primary: 06B23, 06B15;
Secondary: 06A15, 06B05.

\medskip

\textit{Key words and phrases}: compatible tolerance, Dedekind-Mac Neille
completion of a lattice{\small , }weak ordered relation, formal context,
concept lattice.

\bigskip

\noindent \textsc{1. Introduction}

\bigskip

A binary relation on a (complete) lattice $L$ is called \emph{(completely)
compatible} if it is a (complete) sublattice of the direct product $%
L^{2}=L\times L$. A reflexive symmetric and (completely) compatible relation
$T\subseteq L^{2}$ is a called a \emph{(complete) tolerance} of $L$. All
tolerances of a lattice $L$, denoted by Tol$(L)$ form an algebraic lattice
(with respect to the inclusion).

Let $T\in \ $Tol$(L)$ and $X\subseteq L$, $X\neq \emptyset $. If $%
X^{2}\subseteq T$, then the set $X$ is called a preblock of $T$. Blocks are
maximal preblocks (with respect to $\subseteq $). It is known that the
blocks of any tolerance $T$ are convex sublattices of $L$. In [Cz] G. Cz\'{e}%
dli proved that the blocks of $T$ can be ordered in such a way that they
form a lattice. This lattice is denoted by $L/T$ and it is called \emph{the}
\emph{factor lattice of }$L$\emph{\ modulo }$T$\emph{.} The notion of factor
lattices constructed with his method constitute a natural generalization of
that of factor lattices by congruences.

\bigskip

\noindent \textbf{Definition 1.1. }\textit{We say that a binary relation }$R$
\textit{is a weak ordered relation on the lattice }$L$\textit{\ if it
satisfies the following conditions:}

\begin{enumerate}
\item[(1)] \textit{for any }$u,x,y,z\in L$\textit{, }$u\leq x$\textit{, }$%
(x,y)\in R$\textit{\ and }$y\leq z$\textit{\ imply }$(u,z)\in R$\textit{;}

\item[(2)] \textit{given any }$t\in L$\textit{\ and any nonempty finite }$%
A\subseteq L$\textit{, if }$(a,t)\in R$\textit{\ holds for each }$a\in A$%
\textit{\ then }$\left( \bigvee A,t\right) \in R$\textit{;}

\item[(3)] \textit{given any }$z\in L$\textit{\ and any nonempty finite }$%
A\subseteq L$\textit{, if }$(z,a)\in R$\textit{\ holds for each }$a\in A$%
\textit{\ then }$\left( z,\bigwedge A\right) \in R$\textit{.}
\end{enumerate}

\noindent An \emph{ordered relation }$R$ on a complete lattice $L$ is a weak
ordered relation which satisfies conditions (2) and (3) for arbitrary (i.e.
even infinite or empty) $A\subseteq L$. This notion was introduced by S.
Valentini [V], who pointed out that any ordered relation is a completely
compatible relation on $L$. We will show that reflexive weak ordered
relations of $L$ can be characterized as compatible reflexive relations $%
R\subseteq L^{2}$ satisfying $R=\ \leqo\circ R\circ \leqo$ (see [KR1]).
Moreover, we will see that for any $T\in \ $Tol$(L)$, $R:=\ \leqo\circ
T\circ \leqo$ is a weak ordered relation with the property that $T=R\cap
R^{-1}$, where $R^{-1}$ stands for the \emph{inverse relation} of $R$. The
set of weak ordered relations and that of reflexive weak ordered relations
of a lattice $L$ will be denoted by WOR$(L)$ and ReWOR$(L)$, respectively.

The main results of the paper point out the connection between the weak
ordered relations and factor lattices defined by tolerances. It is proved
that for any tolerance $T$ of a lattice $L$ the Dedekind-Mac Neille
completion of $L/T$ is isomorphic to the concept lattice $\mathcal{L}(L,L,R)$
of the context $(L,L,R)$, where $R:=\leqo\circ T\circ \leqo$. It is also
shown that the blocks of $T$ correspond exactly to the concepts $(A,B)\in
\mathcal{L}(L,L,R)$ having the property that $A\cap B\neq \emptyset $. This
result generalizes a result of [KR2], where for any complete lattice $L$ and
any complete tolerance $T\subseteq L^{2}$ the isomorphism $L/T\cong \mathcal{%
L}(L,L,\leqo\circ T\circ \leqo)$ was established.

The paper is structured as follows: In Section 2 some basic notions and the
interrelation between the lattices Tol$(L)$ and ReWOR$(L)$ are presented; In
Section 3, the concept lattice $\mathcal{L}(L,L,\leqo\circ T\circ \leqo)$ is
described and the main results of the paper are presented.

\bigskip

\noindent \textsc{2. Reflexive weak ordered relations and FCA notions}

\bigskip

First, observe that weak ordered relations on a lattice $L$ are also
compatible relations. Indeed, let $R\in \ $WOR$(L)$ and $%
(x_{1},y_{1}),(x_{2},y_{2})\in R$. Condition (1) of Definition 1.1 implies $%
(x_{1},y_{1}\vee y_{2}),(x_{2},y_{1}\vee y_{2})\in R$ and $(x_{1}\wedge
x_{2},y_{1}),(x_{1}\wedge x_{2},y_{2})\in R$. Now, by using conditions (2)
and (3) we obtain $(x_{1}\vee x_{2},y_{1}\vee y_{2})\in R$ and $(x_{1}\wedge
x_{2},y_{1}\wedge y_{2})\in R$, proving the compatibility of $R$.

The relations $\leq $ and $\triangledown =L\times L$ are examples of
reflexive ordered relations. We consider the empty relation $\emptyset $
also a weak ordered relation. Clearly, relation $\emptyset $ is not
reflexive. An other weak ordered relation which is not reflexive in general,
is given in the following

\bigskip

\noindent \textbf{Example} \textbf{2.1. }\textit{Let }$f\colon L\rightarrow
L $\textit{\ be a join-endomorphism of the lattice }$L$\textit{. Then the
binary relation }%
\begin{equation*}
R^{f}:=\{(x,y)\in L^{2}\mid f(x)\leq y\}
\end{equation*}%
\textit{is a weak ordered relation on the lattice }$L$\textit{.}

\bigskip

Indeed, condition (1) of definition 1.1 is satisfied obviously. In order to
prove condition (2), let $A=\{a_{1},...a_{n}\}\subseteq L$ , $t\in L$ and
suppose that $(a_{i},t)\in R^{f}$, for all $i=1,...,n$. Then $f(a_{i})\leq t$%
, $i=1,...,n$ yields $f\left( \underset{i=1}{\overset{n}{\bigvee }}%
a_{i}\right) =\underset{i=1}{\overset{n}{\bigvee }}f(a_{i})\leq t$, proving $%
\left( \bigvee A,t\right) \in R^{f}$; condition (3) is proved analogously.

\bigskip

\noindent \textbf{Proposition 2.2.} \emph{Let }$L$\emph{\ be a lattice. Then}
$($WOR$(L),\subseteq )$ \emph{is a complete lattice, and} $($WOR$(L),\circ )$
i\emph{s a monoid with unit element }$\leq $\emph{. In addition the relations%
}%
\begin{equation}
(R_{1}\cap R_{2})\circ S=\left( R_{1}\circ S\right) \cap \left( R_{2}\circ
S\right)  \tag{D1}
\end{equation}%
\begin{equation}
S\circ (R_{1}\cap R_{2})=\left( S\circ R_{1}\right) \cap \left( S\circ
R_{2}\right)  \tag{D2}
\end{equation}

\noindent \emph{hold for any }$R_{1},R_{2},S\in \ $WOR$(L)$\emph{, that is, }%
$\circ $ \emph{distributes over intersections from both sides.}

\bigskip

\noindent \textit{Proof.} Let $R_{i}\in \ $WOR$(L)$, $i\in I$. It is easy to
check that $\underset{i\in I}{\bigcap }R_{i}$ satisfies the conditions of
Definition 1.1. Indeed, for any $u,x,y,z\in L$\textit{, }$u\leq x$\textit{, }%
$(x,y)\in \ \underset{i\in I}{\bigcap }R_{i}$ and $y\leq z$ imply $u\leq x$%
\textit{, }$(x,y)\in R_{i}$, $y\leq z$, and condition (1) yields $(u,z)\in
R_{i}$, $i\in I$. Hence we get $(u,z)\in \ \underset{i\in I}{\bigcap }R_{i}$%
, i.e. $\underset{i\in I}{\bigcap }R_{i}$ satisfies condition (1).

Now let $t\in L$ and $A\subseteq L$ a finite nonempty set. In order to check
condition (2), assume that $(a,t)\in \ \underset{i\in I}{\bigcap }R_{i}$,
for each $a\in A$. Then $(a,t)\in R_{i}$, for all $i\in I$ and each $a\in A$%
, and hence condition (2) of Definition 1.1 yields $\left( \bigvee
A,t\right) \in R_{i}$, for all $i\in I$. Thus we get $\left( \bigvee
A,t\right) \in \ \underset{i\in I}{\bigcap }R_{i}$, and this mean that
condition (2) holds for $\underset{i\in I}{\bigcap }R_{i}$ also. The fact
that $\underset{i\in I}{\bigcap }R_{i}$ satisfies condition (3) is proved
similarly.

Thus we have proved that $\underset{i\in I}{\bigcap }R_{i}\in \ $WOR$(L)$.
Since the relation $\nabla $ is the greatest element in WOR$(L)$, we obtain
that $($WOR$(L),\subseteq )$ is a complete lattice.

Now let $R,S\in \ $WOR$(L)$. We prove that $R\circ S\in \ $WOR$(L)$. Indeed,
condition (1) of Definition 1.1 is satisfied trivially. In order to prove
condition (2), take any finite $A\subseteq L$, $A\neq \emptyset $ and $t\in
L $ and assume that $(a,t)\in R\circ S$, for each $a\in A$. Then for each $%
a\in A$ there exists a $z_{a}\in L$ such that $(a,z_{a})\in R$ and $%
(z_{a},t)\in S$. Since $R$ is a compatible relation, we$S$ obtain $(\bigvee
A,\bigvee z_{a})\in R$. Applying condition (2) for $S$ we get $\left(
\bigvee z_{a},t\right) \in S$. Thus we deduce $\left( \bigvee A,t\right) \in
R\circ S$, proving that relation $R\circ S$ satisfies condition (2). The
fact that $R\circ S$ satisfies condition (3) is proved dually. Hence $R\circ
S\in \ $WOR$(L)$. Because $\circ $ is an associative operation, $($WOR$%
(L),\circ )$ is a semigroup.

We already noted that $\leq $ belongs to WOR$(L)$. By using condition (1),
we obtain that $\leqo\circ R\subseteq R$ and $R\circ \leqo\subseteq R$ holds
for any $R\in \ $WOR$(L)$. Since $\leq $ is a reflexive relation, the
inclusions $R\subseteq \leqo\circ R$ and $R\subseteq R\circ \leqo$ are
obvious. Thus we have%
\begin{equation}
\leqo\circ R=R\circ \leqo\ =R\text{,}  \tag{U}
\end{equation}

\noindent and this means that $($WOR$(L),\circ )$ is a monoid with unit
element $\leq $.

Next, take any $R_{1},R_{2},S\in \ $WOR$(L)$, and prove\emph{\ }identity
(D1). The inclusion $(R_{1}\cap R_{2})\circ S\subseteq \left( R_{1}\circ
S\right) \cap \left( R_{2}\circ S\right) $ is obvious. In order to prove the
converse inclusion, take any $(x,y)\in \left( R_{1}\circ S\right) \cap
\left( R_{2}\circ S\right) $. Then there exist some $z_{1},z_{2}\in L$ such
that $(x,z_{1})\in R_{1}$, $(x,z_{2})\in R_{2}$ and $(z_{1},y),(z_{2},y)\in
S $. Then, in view of condition (1), $z_{1},z_{2}\leq z_{1}\vee z_{2}$
implies $(x,z_{1}\vee z_{2})\in R_{1}\cap R_{2}$, and applying condition (2)
for $S$ we get $(z_{1}\vee z_{2},y)\in S$. Hence we obtain $(x,y)\in
(R_{1}\cap R_{2})\circ S$, proving identity (D1). Identity (D2) is proved
similarly.\hfill $\square $

\bigskip

As an immediate consequence of Proposition 2.2 we obtain

\bigskip

\noindent \textbf{Corollary 2.3.} \emph{Let }$R$\emph{\ be a binary relation
on the lattice }$L$\emph{. Then the following are equivalent:}

\begin{enumerate}
\item[(i)] $R$\emph{\ is a reflexive weak ordered relation;}

\item[(ii)] $R$\emph{\ is a reflexive compatible relation on }$L$\emph{\
which satisfies }$\leqo\circ R\circ \leqo=R$\emph{.}
\end{enumerate}

\bigskip

\noindent \textit{Proof.} (i)$\Rightarrow $(ii). We have already shown that
any weak ordered relation $R$ is compatible. Since $\leq $ is the unit of
the monoid $($WOR$(L),\circ )$, $\leqo\circ R\circ \leqo=R$\emph{\ }is clear.

\noindent (ii)$\Rightarrow $(i). Because $R$ is a reflexive compatible
relation, for any $z,t\in L$ and any finite $A=\{a_{1},...,a_{n}\}\subseteq
L $, $(a_{i},t)$, $(z,a_{i})\in R$, for all $i=1,...,n$ imply $(a_{1}\vee
...\vee a_{n},t\vee ...\vee t)\in R$ and $(z\wedge ...\wedge z,a_{1}\wedge
...\wedge a_{n})\in R$. Thus \textit{\ }$\left( \bigvee A,t\right) \in R$%
\textit{\ }and $\left( z,\bigwedge A\right) \in R$ hold, proving that
conditions (2) and (3) are satisfied by $R$. In order to prove condition (1)
take any $u,x,y,z\in L$\textit{, }with $u\leq x$\textit{, }$(x,y)\in R$%
\textit{\ and }$y\leq z$. Then $(u,z)\in \leqo\circ R\circ \leqo=R$, and
this proves condition (1). Hence $R$ is a reflexive weak ordered relation.
\hfill $\square $

\bigskip

Let us denote the set of compatible reflexive relations of a lattice $L$ by
Re$(L)$. It was proved in [PR] that Re$(L)$ forms an algebraic lattice with
respect to $\subseteq $. Clearly, the least element of Re$(L)$ is the
identity relation on $L$, i.e. $\triangle =\{(x,x)\mid x\in L\}$. The
following lemmas contain some properties of Re$(L)$ and Tol$(L)$ which will
be useful in our proofs.

\bigskip

\noindent \textbf{Lemma 2.4.} \emph{Let }$L$\emph{\ be a lattice and }$%
R_{1},R_{2},S\in \ $Re$(L)$. \emph{Then the following assertions hold true:}

\begin{enumerate}
\item[(a)] $\left( R_{1}\circ R_{2}\right) \cap S\subseteq (R_{1}\cap
S)\circ (R_{2}\cap S)$.

\item[(b)] $(R_{1}\cap R_{2})\circ S=\left( R_{1}\circ S\right) \cap \left(
R_{2}\circ S\right) $ \emph{and} $S\circ (R_{1}\cap R_{2})=\left( S\circ
R\right) _{1}\cap \left( S\circ R\right) _{2}$.

\item[(c)] \emph{For arbitrary} $T_{1},T_{2}\in \ $Tol$(L)$, \emph{we have} $%
T_{1}=T_{2}\Leftrightarrow T_{1}\cap \leqo=\ T_{2}\cap \leqo$.
\end{enumerate}

\bigskip

\noindent We note that relation (a) above is proved in [ChR], while (b) and
(c) can be found in [KR1], moreover, both (a) and (b) are valid even in any
algebra with a (ternary) majority term.

\bigskip

\noindent \textbf{Lemma 2.5.}\emph{\ Let} \emph{be a lattice, }$S\in \ $ReWOR%
$(L)$ and $T\in \ $Tol$(L)$\emph{. Then}

\begin{enumerate}
\item[(i)] $S$ includes $\leq $ and $\leqo\circ S^{-1}=S^{-1}\circ \leqo%
=\triangledown $;

\item[(ii)] $\leqo\circ T\circ \leqo\in \ $ReWOR$(L)$;

\item[(iii)] $T=\left( \leqo\circ T\circ \leqo\right) \cap \left( \geqo %
\circ T\circ \geqo \right) $.
\end{enumerate}

\bigskip

\noindent \textit{Proof.} (i) Since $S$ is reflexive, $\leqo\subseteq S$ is
clear. Then $\geqo\subseteq S^{-1}$ also holds. Now take an arbitrary $%
(x,y)\in L^{2}$. Then $x\leq x\vee y$, and $x\vee y\geq y$ yields $(x\vee
y,y)\in S^{-1}$. Hence $(x,y)\in \leqo\circ S^{-1}$, and this proves $\leqo%
\circ S^{-1}=\triangledown $. Since $\left( x,x\wedge y\right) \in S^{-1}$
and $x\wedge y\leq y$, the relation $S^{-1}\circ \leqo=\triangledown $ holds
also.

(ii) Clearly, $R:=\leqo\circ T\circ \leqo$ is a reflexive compatible
relation on $L$, and $\leqo\circ R\circ \leqo=R$, because\emph{\ }$\leqo%
\circ \leqo\ =\ \leq $. Hence, in view of Corollary 2.3, $R$ is a weak
ordered relation.

(iii) Denote $S:=\left( \leqo\circ T\circ \leqo\right) \cap \left( \geqo%
\circ T\circ \geqo\right) $. Since $\geqo\circ T\circ \geqo=\left( \leqo%
\circ T\circ \leqo\right) ^{-1}$, $S$ is symmetric and hence $S\in \ $Tol$%
(L) $. Now, in view of Lemma 2.4.(c), to prove (ii) it is enough to show
that $T\cap \leqo=\ S\cap \leqo$. Since $T\subseteq S$, the inclusion $T\cap %
\leqo\subseteq \ S\cap \leqo$ is clear. On the other hand, $S\cap \leqo%
=\left( \leqo\circ T\circ \leqo\right) \cap \leqo\cap (\geqo\circ T\circ \ %
\geqo)$. As $\leqo\circ T\circ \leqo$ $\in \ $ReWOR$(L)$, in view of (i) it
includes $\leq $, and hence $\left( \leqo\circ T\circ \leqo\right) \cap \leqo%
=\leqo$. Thus we get $S\cap \leqo=\leqo\cap \left( \geqo\circ T\circ \geqo%
\right) $. By using Lemma 2.4(a) we obtain $\leqo\cap \left( \geqo\circ
T\circ \geqo\right) \subseteq (\leqo\cap \geqo)\circ \left( \leqo\cap
T\right) \circ \ (\leqo\cap \geqo)=\triangle \circ \left( \leqo\cap T\right)
\circ \triangle =\leqo\cap T$. Hence $S\cap \leqo\subseteq T\cap \leqo$, and
this proves $T\cap \leqo=\ S\cap \leqo$. \hfill $\square $

\bigskip

\noindent \textbf{Theorem 2.6.}\emph{\ The mappings }

$\alpha \colon $ReWOR$(L)\rightarrow \ $Tol$(L)$, $\alpha (R)=R\cap R^{-1}$
\emph{and}

$\beta \colon $Tol$(L)\rightarrow \ $ReWOR$(L)$, $\beta (T)=\leqo\circ
T\circ \leqo$

\noindent \emph{are lattice isomorphisms and they are inverses each of other}%
.

\bigskip

\noindent \textit{Proof.} Obviously, if $R$ is a compatible reflexive
relation on the lattice $L$ then $\alpha (R)=R\cap R^{-1}\in \ $Tol$(L)$. In
view of Lemma 2.5.(ii), for any $T\in \ $Tol$(L)$, we have $\beta (T)=\leqo%
\circ T\circ \leqo\in \ $ReWOR$(L)$. Thus the maps $\alpha $ and $\beta $
are correctly defined. First, we prove that they are inverses each of other:

Indeed, by Lemma 2.4(iii) we get

\noindent $\alpha \left( \beta (T)\right) =\left( \leqo\circ T\circ \leqo%
\right) \cap \left( \leqo\circ T\circ \leqo\right) ^{-1}=\left( \leqo\circ
T\circ \leqo\right) \cap \left( \geqo\circ T\circ \geqo\right) =T$, for each
$T\in $Tol$(L)$.

On the other hand, we have $\beta (\alpha (R))=\leqo\circ \left( R\cap
R^{-1}\right) \circ \leqo$, for any $R\in \ $ReWOR$(L)$. We are going to
prove $\beta (\alpha (R))=R$, i.e. $\leqo\circ \left( R\cap R^{-1}\right)
\circ \leqo=R$.

\noindent As $\leq $, $R$ and $R^{-1}$ belong to Re$(L)$, by applying Lemma
2.4(b) we obtain:

$\leqo\circ \left( R\cap R^{-1}\right) \circ \leqo=\left( \left( \leqo\circ
R\right) \cap \left( \leqo\circ R^{-1}\right) \right) \circ \leqo$.

\noindent Since $R\in \ $ReWOR$(L)$, we get $\leqo\circ R=R$ and $\leqo\circ
R^{-1}=\nabla $, according to Proposition 2.2 and Lemma 2.5(i). Summarizing
we obtain:

$\left( \left( \leqo\circ R\right) \cap \left( \leqo\circ R^{-1}\right)
\right) \circ \leqo\ =(R\cap \nabla )\circ \leqo\ =R\circ \leqo=R$,

\noindent because $R\circ \leqo=R$ is also true. Thus $\beta (\alpha (R))=R$%
, and hence $\beta =\alpha ^{-1}$.

Finally, observe that both $\alpha $ and $\beta $ are order-preserving,
because so are $\circ $, $\cap $, and taking inverses. Thus $\alpha $ and $%
\beta $ are lattice isomorphisms. \hfill $\square $

\bigskip

The following corollary is obvious:

\bigskip

\noindent \textbf{Corollary 2.7.}\emph{\ Any reflexive weak ordered relation
}$R\subseteq L^{2}$ \emph{can be represented in the form }$R=\leqo\circ
\left( R\cap R^{-1}\right) \circ \leqo$\emph{, in other words, it can be
derived from a compatible tolerance }$T=R\cap R^{-1}$\emph{\ of the lattice }%
$L$\emph{.}

\bigskip

\textbf{Some notions from Formal Concept Analysis}

\bigskip

A \emph{formal context} is a triple $\mathcal{K=}(G,M,I)$, where $G$ and $M$
are sets and $I\subseteq G\times M$ is a binary relation. The basic notions
of Formal Concept Analysis (FCA) can be found e.g. in [GW] or [W]. By
defining for all subsets $A\subseteq G$ and $B\subseteq M$

\begin{center}
$A^{I}:=\{m\in M\mid (g,m)\in I$, for all $g\in A\}$,

$^{I}B:=\{g\in G\mid (g,m)\in I$, for all $m\in B\}$
\end{center}

\noindent we establish a Galois connection between the power-set lattices $%
\wp (G)$ and $\wp (M)$. We will use the notations $^{I}(A^{I})=A^{II}$ and $%
\left( ^{I}B\right) ^{I}=B^{II}$, for any $A\subseteq G$ and $B\subseteq M$.
The obtained maps $A\rightarrow A^{II}$, $A\subseteq G$ \ and $B\rightarrow
B^{II}$, $B\subseteq M$ are closure operators on $\wp (G),$ respectively $%
\wp (M)$.

A \emph{formal concept} of the context $\mathcal{K}$ is a pair $(A,B)\in \wp
(G)\times \wp (M)$ with $A^{I}=B$ and $^{I}B=A$, where the set $A$ is called
the \emph{extent} and $B$ is called the \emph{intent} of the concept $(A,B)$%
. It is easy to check that a pair $(A,B)\in \wp (G)\times \wp (M)$ is a
concept if and only if $(A,B)=(A^{II},A^{I})=(^{I}B,B^{II})$. The concepts
of the context $(G,M,I)$ can be also characterized as those pairs $(A,B)\in
\wp (G)\times \wp (M)$ whose products are maximal sets with the property $%
A\times B\subseteq I$. The set of all concepts of the context $\mathcal{K}$
is denoted by $\mathcal{L}(\mathcal{K)}$. This set $\mathcal{L}(G,M,I%
\mathcal{)}$ ordered by the relation $\leq $ defined as follows

\begin{equation*}
(A_{1},B_{1})\leq (A_{2},B_{2})\Leftrightarrow A_{1}\subseteq
A_{2}\Leftrightarrow B_{1}\supseteq B_{2}\text{,}
\end{equation*}

\noindent forms a complete lattice, called the \emph{concept lattice of the
context} $\mathcal{K}=(G,M,I)$.

Let us consider the concepts $\gamma (x)=(\{x\}^{II},\{x\}^{I})$ and $\mu
(y)=(^{I}\{y\},\{y\}^{II})$, for any $x\in G$ and $y\in M$. It can be easily
proved that for any concept $(A,B)\in \mathcal{L}(G,M,I\mathcal{)}$, we have
in $\mathcal{L}(G,M,I\mathcal{)}$:%
\begin{equation}
(A,B)=\bigvee \{\gamma (x)\mid x\in A\}=\bigwedge \{\mu (y)\mid y\in B\}%
\text{.}  \tag{E}
\end{equation}

\noindent The following assertion is a part of Basic Theorem on Concept
Lattices from [GW; Thm. 3.]:

\bigskip

\noindent \textbf{Proposition 2.8.} \emph{A}\textbf{\ }\emph{complete
lattice }$L$ \emph{is isomorphic to }$\mathcal{L}(G,M,I\mathcal{)}$ \emph{if
and only if there are some mappings }$\widetilde{\gamma }\colon G\rightarrow
L$\emph{\ and }$\widetilde{\mu }\colon M\rightarrow L$ \emph{such that}
\emph{the set} $\{\widetilde{\gamma }(g)\mid g\in G\}$ \emph{is
supremum-dense in} $L$, $\{\widetilde{\mu }(m)\mid m\in M\}$ \emph{is
infimum-dense in} $L$ \emph{and }$(g,m)\in I$ \emph{is equivalent to} $%
\widetilde{\gamma }(g)\leq \widetilde{\mu }(m)$ \emph{for all }$g\in G$\emph{%
\ and }$m\in M$\emph{.}

\bigskip

\noindent \textsc{3. Concept lattices induced by weak ordered relations}

\bigskip

\noindent \textbf{Proposition 3.1.}\emph{\ Let }$L$ \emph{be a lattice, }$%
R\subseteq L^{2}$ \emph{a weak ordered relation, and }$(A,B)$\emph{\ a
concept of the context }$(L,L,R)$\emph{. Then }$A$\emph{\ is an ideal and }$%
B $\emph{\ is a filter in} $L$.

\bigskip

\noindent \textit{Proof. }Suppose that $x\leq a$ for some $a\in A$ and $x\in
L$. Since $(a,b)\in R$ for all $a\in A$ and $b\in B$ and $R$ is a weak
ordered relation, we obtain $(x,b)\in B$, for all $b\in B$. Hence $x\in \
^{I}B=A$. Now let $a_{1},a_{2}\in A$. Then for each $b\in B$ the relations $%
\left( a_{1},b\right) ,\left( a_{2},b\right) \in R$ imply $(a_{1}\vee
a_{2},b)\in R$ (see Definition 1.1(2)). Hence $a_{1}\vee a_{2}\in \ ^{I}B=A$%
. This proves that $A$ is an ideal of $L$. The fact that $B$ is a filter of $%
L$ is proved dually. \hfill $\square $

\bigskip

For any subset $X\subseteq L$ of a lattice $L$, let $[X)$ and $(X]$ denote
the filter and the ideal generated by $X$, respectively. We will use the
following (see also [Cz; Lemma 2]):

\bigskip

\noindent \textbf{Lemma 3.2.} (Gr\"{a}tzer [G; Lemma 6]) \emph{For any
convex sublattice }$C$\emph{\ of the lattice }$L$\emph{\ the equality} $%
C=(C]\cap \lbrack C)$ \emph{holds. Moreover, if }$C$\emph{\ is the
intersection of an ideal }$I$\emph{\ and of a filter }$F$\emph{\ of }$L$%
\emph{, then} $I=(C]$ \emph{and} $F=[C)$.

\bigskip

\noindent \textbf{Proposition 3.3.} \emph{Let }$T$ \emph{be a tolerance of
the lattice }$L$. \emph{Then for any block }$C$ \emph{of} $T$ \emph{the pair}
$\left( (C],[C)\right) $\emph{\ coincides to the unique concept} $(A,B)\in
\mathcal{L}(L,L,\leqo\circ T\circ \leqo)$ \emph{with} $C=A\cap B$. \emph{For
any} \emph{concept }$(A,B)\in \mathcal{L}(L,L,\leqo\circ T\circ \leqo)$, $%
A\cap B$ \emph{is a block of} $T$ \emph{whenever} $A\cap B\neq \emptyset $.

\bigskip

\noindent \textit{Proof.} Denote $R:=\leqo\circ T\circ \leqo$. Then, in view
of Theorem 2.6, $T=R\cap R^{-1}$. Let $C$ be a block of $T$. Then $C\times
C\subseteq T\subseteq R$, hence $C\subseteq C^{R}$. Clearly, $%
(C^{RR},C^{R})\in \mathcal{L}(L,L,R)$ and hence $C^{RR}\times C^{R}\subseteq
R$. Since $C\subseteq C^{RR}$ always holds, we obtain $C\subseteq C^{RR}\cap
C^{R}$. We claim that $C=C^{RR}\cap C^{R}$. Indeed, $\left( C^{RR}\cap
C^{R}\right) \times \left( C^{RR}\cap C^{R}\right) \subseteq \left(
C^{RR}\times C^{R}\right) \cap (C^{R}\times C^{RR})=R\cap R^{-1}=T$, and
this means that $C^{RR}\cap C^{R}$ is a preblock of $T$. Since $C$ is block
and $C\subseteq C^{RR}\cap C^{R}$, we obtain $C=C^{RR}\cap C^{R}$. Because $%
R $ is a weak ordered relation and $(C^{RR},C^{R})$ is a concept of the
context $(L,L,R)$, the extent $C^{RR}$ is an ideal of $L$ and the intent $%
C^{R}$ is a filter of $L$, according to Proposition 3.1. As $C$ is a convex
sublattice of $L$, by using Lemma 3.2, we obtain $C^{RR}=(C]$, $C^{R}=[C)$
and $C=(C]\cap \lbrack C)$. Then $\left( (C],[C)\right) \in \mathcal{L}(L,L,%
\leqo\circ T\circ \leqo)$ also holds, so $C$ is of the required form. Now,
assume that $C=A\cap B$ holds for some concept $(A,B)\in \mathcal{L}(L,L,R)$%
. Since $A$ is an ideal and $B$ is a filter of $L$, in view of Lemma 3.2. we
obtain $(A,B)=\left( (C],[C)\right) $. Therefore $\left( (C],[C)\right) $ is
the unique concept with the required property.

Finally, take any $(A,B)\in \mathcal{L}(L,L,R)$ such that $D:=A\cap B\neq
\emptyset $. Since $A$ is an ideal and $B$ is a filter of $L$ according to
Proposition 3.1, we get that $D$ is a convex sublattice of $L$. Then in view
of Lemma 3.2, we have $A=(D]$ and $B=[D)$. Since $D\times D\subseteq \left(
A\times B\right) \cap (B\times A)\subseteq R\cap R^{-1}=T$, $D$ is preblock
of $T$. Then there exists (at least one) block $C$ of $T$ such that $%
D\subseteq C$. Then $\left( (C],[C)\right) $ is a concept of the context $%
(L,L,R)$, moreover, we have $A=(D]\subseteq (C]$ and $B=[D)\subseteq \lbrack
C)$. Since $(A,B)$ is also a concept of the same context, these relations
imply $(A,B)=\left( (C],[C)\right) $. Then $D=A\cap B=(C]\cap \lbrack C)=C$.
This proves that $A\cap B$ is a block of $T$. \hfill $\square $

\bigskip \pagebreak

\noindent \textbf{Proposition 3.4. }\emph{Let }$T$ \emph{be a tolerance on
the lattice }$L$. \emph{Then the mapping }%
\begin{equation*}
\delta \colon L/T\rightarrow \mathcal{L}(L,L,\leqo\circ T\circ \leqo)\text{,
}\delta (C)=\left( (C],[C)\right) \text{, }C\in L/T
\end{equation*}

\noindent \emph{is a lattice embedding.}

\bigskip

\noindent P\textit{roof. }In view of Proposition 3.3, for any tolerance
block $C\in L/T$, we have $\left( (C],[C)\right) \in \mathcal{L}(L,L,\leqo%
\circ T\circ \leqo)$, i.e. the mapping $\delta $ is well-defined. Assume
that $\delta (C)=\delta (D)$ for some $C,D\in L/T$. Then $\left(
(C],[C)\right) =\left( (D],[D)\right) $ and Proposition 3.3 imply $C=(C]\cap
\lbrack C)=(D]\cap \lbrack D)=D$, proving that $\delta $ is one to one.

Further, denote $R:=\leqo\circ T\circ \leqo$, and assume that for some $%
B_{1},B_{2},E,F\in L/T$ the equalities $B_{1}\vee B_{2}=E$ and $B_{1}\wedge
B_{2}=F$ hold in the factor lattice $L/T$. Now, in view of [Cz; Lemma 4] we
have:

$(B_{1}\cup B_{2}]\subseteq (E]$, $[B_{1})\cap \lbrack B_{2})=[E)$ and

$(B_{1}]\cap (B_{2}]=(F]$, $[B_{1}\vee B_{2})\subseteq \lbrack F)$.

\noindent Then $\delta (B_{1})=\left( (B_{1}],[B_{1})\right) $, $\delta
(B_{2})=\left( (B_{2}],[B_{2})\right) $ and $\delta (B_{1}\vee B_{2})=\left(
(E],[E)\right) =\left( (^{R}[E),[E)\right) $, $\delta (B_{1}\wedge
B_{2})=\left( (F],[F)\right) =\left( (F],(F]^{R}\right) $. On the other
hand, in view of [GW] (or [W]), the $\vee $ and $\wedge $ operation in $%
\mathcal{L}(L,L,R)$ have the form

$\left( (B_{1}],[B_{1})\right) \vee \left( (B_{2}],[B_{2})\right) =\left(
^{R}([B_{1})\cap \lbrack B_{2})\right) ,[B_{1})\cap \lbrack B_{2}))$

\noindent and

$\left( (B_{1}],[B_{1})\right) \wedge \left( (B_{2}],[B_{2})\right)
=((B_{1}]\cap (B_{2}],((B_{1}]\cap (B_{2}])^{R})$.

\noindent Since $[E)=[B_{1})\cap \lbrack B_{2})$, we obtain

$\delta (B_{1}\vee B_{2})=\left( (^{R}[E),[E)\right) =\left(
(B_{1}],[B_{1})\right) \vee \left( (B_{2}],[B_{2})\right) =\delta
(B_{1})\vee \delta (B_{2})$,

\noindent and similarly, $(B_{1}]\cap (B_{2}]=(F]$ implies

$\delta (B_{1}\wedge B_{2})=\left( (F],(F]^{R}\right) =\left(
(B_{1}],[B_{1})\right) \wedge \left( (B_{2}],[B_{2})\right) =\delta
(B_{1})\wedge \delta (B_{2})$.

\noindent Thus $\delta $ is a lattice embedding. \hfill $\square $

\bigskip

Denote the Dedekind-Mac Neille completion of a lattice $L$ by DM$(L)$. It is
known, that for any lattice $L$, DM$(L)$ is isomorphic to the concept
lattice $\mathcal{L}(L,L,\leq )$.

\bigskip

\noindent \textbf{Theorem 3.5.}\emph{\ Let }$T$ \emph{be a tolerance of the
lattice }$L$ \emph{and denote} $R:=\leqo\circ T\circ \leqo$. \emph{Then }DM$%
(L/T)$ \emph{is isomorphic to the concept lattice} $\mathcal{L}(L,L,R)$.

\bigskip

\noindent P\textit{roof. }First, observe that for any\textit{\ }$x\in L$, we
have $(x,x)\in R$ and this implies $x\in \{x\}^{R}$ and $x\in \ ^{R}\{x\}$.
As $x\in \{x\}^{RR}$ always holds, we get $\{x\}^{RR}\cap \{x\}^{R}\neq
\emptyset $, and in view of Proposition 3.3 this means that the concept $%
\gamma (x)=\left( \{x\}^{RR},\{x\}^{R}\right) \in $ $\mathcal{L}(L,L,R)$ has
the form $\gamma (x)=((C],[C))$, where $C$ is a block of $T$ such that $%
C=\{x\}^{RR}\cap \{x\}^{R}$. Thus $C$ contains $x$. Similarly is proved that
$^{R}\{x\}\cap \{x\}^{RR}\neq \emptyset $ implies that the concept $\mu
(x)=\left( ^{R}\{x\},\{x\}^{RR}\right) $ has the form $\mu (x)=((D],[D))$,
where $D$ is a block of $T$ such that $D=\ ^{R}\{x\}\cap \{x\}^{RR}\ni x$.
Since, in view of (E), the set is $\{\gamma (x)\mid x\in L\}$ is supremum
dense and the set $\{\mu (x)\mid x\in L\}$ is infimum dense in $\mathcal{L}%
(L,L,R)$, we obtain that the concepts $((B],[B))$, $B\in L/T$ form a set
which is both supremum- and infimum-dense in $\mathcal{L}(L,L,R)$.

Now, consider the lattice $\mathcal{L}(L/T,L/T,\leq )$. Since DM$(L/T)$ is
isomorphic to $\mathcal{L}(L/T,L/T,\leq )$, to prove our theorem it is
enough to show that $\mathcal{L}(L/T,L/T,\leq )\cong \mathcal{L}(L,L,R)$. In
order to apply Proposition 2.8, for any $B\in L/T$ define the mappings $%
\widetilde{\gamma }\colon L/T\rightarrow \mathcal{L}(L,L,R)$, $\widetilde{%
\mu }\colon L/T\rightarrow \mathcal{L}(L,L,R)$ to be equal to the mapping $%
\delta \colon L/T\rightarrow \mathcal{L}(L,L,R),\delta (B)=\left(
(B],[B)\right) $, $B\in L/T$, i.e. let%
\begin{equation*}
\widetilde{\gamma }=\widetilde{\mu }=\delta \text{.}
\end{equation*}

\noindent Then $\{\widetilde{\gamma }(B)\mid B\in L/T\}=\{((B],[B))\mid B\in
L/T\}$ is supremum dense in $\mathcal{L}(L,L,R)$, and $\{\widetilde{\mu }%
(B)\mid B\in L/T\}=\{((B],[B))\mid B\in L/T\}$ is infimum-dense in $V$.

Now suppose that $B\leq C$ holds in $L/T$ for some blocks $B,C\in L/T$. This
is equivalent to $B=B\wedge C$ in $L/T$. Since in view of Proposition 3.4, $%
\delta $ is a lattice embedding, $B=B\wedge C\Leftrightarrow \delta
(B)=\delta (B\wedge C)=\delta (B)\wedge \delta (C)$ in $\mathcal{L}(L,L,R)$.
Since the latter relation is equivalent to $\delta (B)\leq \delta (C)$, we
obtain that $B\leq C$ if and only if $\delta (B)\leq \delta (C)$. Finally,
by applying Proposition 2.8, we obtain $\mathcal{L}(L/T,L/T,\leq )\cong
\mathcal{L}(L,L,R)$ and this completes our proof. \hfill $\square $

\bigskip

\noindent If $L/T$ is a complete lattice, then obviously DM$(L/T)=L/T$.
Hence we obtain:

\bigskip

\noindent \textbf{Corollary 3.6.} \emph{Let }$T$ \emph{be a tolerance of the
lattice }$L$ \emph{such that the factor lattice }$L/T$ \emph{is complete.
Then }$L/T\cong \mathcal{L}(L,L,\leqo\circ T\circ \leqo)$.

\bigskip

\noindent \textbf{Remark 3.7.} \emph{This is the case when the factor
lattice }$L/T$\emph{\ is finite. The same result, i.e. }$L/T\cong L(L,L,\leqo%
\circ T\circ \leqo)$\emph{\ we obtain also for a complete tolerance }$T$%
\emph{\ of a complete lattice }$L$\emph{, because then }$L/T$\emph{\ is a
complete lattice. This isomorphism for complete tolerances is also
established in }[KR2]\emph{. We note that in this case in view of }[KR1]%
\emph{\ and }[KR2]\emph{\ the product }$\leqo\circ T\circ \leqo$\emph{\ is
an ordered relation.}

\bigskip

\noindent \textbf{References}

\bigskip

\noindent \lbrack Cz] Cz\'{e}dli, G.: Factor lattices by tolerances, \textit{%
Acta Sci. Math. (Szeged)}, \textbf{44} (1982), 35--42.

\medskip

\noindent \lbrack ChR] Chajda I., Radeleczki, S.: 0-conditions and tolerance
schemes. \textit{Acta Math. Univ. Comeniane}. N.S. \textbf{72} (2003),
177-184.

\medskip

\noindent \lbrack G] Gr\"{a}tzer, G.: General Lattice Theory,
Akademi-Verlag, Berlin, (1978).

\medskip

\noindent \lbrack GW] Ganter, B. and Wille, R.: Formal Concept Analysis:
Mathematical Foundations, Springer, Berlin/Heidelberg, (1999).

\medskip

\noindent \lbrack KR1] Kaarli, K. and Radeleczki S.: Representation of
integral quantales by tolerances, \textit{Algebra Universalis}, \textbf{79}
(5) (2018), https://doi.org/10.1007/s00012-018-0484-1.

\medskip

\noindent \lbrack KR2] Kaarli, K. and Radeleczki S.: Ordered relations,
concept lattices and self-bonds, \textit{Journal of Multiple-Valued Logic
and Soft Computing}, submitted

\medskip

\noindent \lbrack PR] P\"{o}schel, R., Radeleczki, S.: Related structures
with involution, \textit{Acta Math. Hungar,} \textbf{123} (1-2), (2009),
169-185.

\medskip

\noindent \lbrack V] Valentini, S.: Representation theorems for quantales.
\textit{Math. Logic Quart.} 40, (1994) 182-190.

\medskip

\noindent \lbrack W] Wille, R.: Restructuring lattice theory: an approach
based on hierarchies of concepts. In: I. Rival editor, Ordered Sets, Reidel,
Dordrecht-Boston (1982), 445-470.

\bigskip

\noindent \textsc{S\'{a}ndor Radeleczki}

Institute of mathematics, University of Miskolc, Hungary

\textit{e-mail:} matradi@uni-miskolc.hu

\bigskip


\end{document}